\documentclass[11pt,a4paper]{article}
\usepackage[margin=28mm]{geometry}
\usepackage{amsmath,amssymb,amsthm}
\usepackage{tikz,float}
\usepackage[hidelinks]{hyperref}
\newcommand{\F}{\mathbb F}
\newcommand{\Z}{\mathbb Z}

\newcommand{\PP}{\mathfrak P}
\newcommand{\Tr}{\operatorname{tr}}

\newtheorem{theorem}{Theorem}
\newtheorem{lemma}[theorem]{Lemma}
\newtheorem{corollary}[theorem]{Corollary}

\newtheorem{result}[theorem]{Result}
\title{Projective planes with abelian collineation groups of type (F) are Desarguesian}
\author{Jamie Jiaqi Low and Bernhard Schmidt\\[1ex]
\small Division of Mathematical Sciences, School of Physical and Mathematical Sciences\\
\small Nanyang Technological University, 21 Nanyang Link, Singapore 637371\\
\small \texttt{M230143@e.ntu.edu.sg}, \texttt{bernhard@ntu.edu.sg}}
\date{}

\begin{document}
\maketitle

\begin{abstract}
We use prime ideal factorization and residue class fields of cyclotomic integer rings to 
prove that every abelian direct product difference set is classical. 
This implies that every projective plane of order $n$ with an abelian 
collineation group of order $n(n-1)$ is Desarguesian, proving 
a conjecture of Jungnickel and de Resmini.
As a  consequence, we obtain the classification of circular Costas maps
and of Costas polynomials conjectured by Muratovi\'c-Ribi\'c, Pott, Thomson and Wang
and by Torres and Rubio.
\end{abstract}



\section{Introduction}

A finite \emph{projective plane} $\Pi=(\mathcal P,\mathcal L)$ consists
of a set $\mathcal P$ of points and a set $\mathcal L$ of lines, regarded
as subsets of $\mathcal P$, such that any two distinct points lie on
exactly one line, any two distinct lines meet in exactly one point,
and there are four points no three of which are collinear.
The plane has \emph{order} $n$ if every line contains $n+1$ points.
Then every point lies on $n+1$ lines and
$|\mathcal P|=|\mathcal L|=n^2+n+1$.

A \emph{collineation} of $\Pi$ is a permutation of $\mathcal P$ that
maps lines to lines. The group of all collineations of $\Pi$, under
composition, is denoted by $\operatorname{Aut}(\Pi)$.
A group $G\leq\operatorname{Aut}(\Pi)$ is
\emph{quasiregular} if its induced permutation group on each point
orbit and each line orbit is regular. Here regular means transitive
with trivial point stabilizers in the induced permutation group.
Write $P^g$ and $\ell^g$ for the images of a point and a line under $g$,
and $P^G=\{P^g:g\in G\}$ and $\ell^G=\{\ell^g:g\in G\}$ for their orbits.
Explicitly,
quasiregularity means that, for every $g\in G$,
\[
\begin{aligned}
P^g=P&\ \Longrightarrow\ Q^g=Q
&&\text{for every }Q\in P^G,\\
\ell^g=\ell&\ \Longrightarrow\ m^g=m
&&\text{for every }m\in \ell^G.
\end{aligned}
\]
Note that every abelian collineation group is quasiregular.

For $n>2$, a quasiregular collineation group $G$ is of \emph{type~(F)}
if $|G|=n(n-1)$, it has five point orbits, and its fixed points and
fixed lines are precisely two distinct points $P,Q$ and two distinct
lines $\ell,m$ with
\[
\ell=PQ,\qquad P\in m,\qquad Q\notin m.
\]
The point orbits are then
\[
\{P\},\quad\{Q\},\quad\ell\setminus\{P,Q\},\quad
m\setminus\{P\},\quad\mathcal P\setminus(\ell\cup m),
\]
of sizes $1,1,n-1,n,n(n-1)$, respectively.
This is the type~(F) configuration of Dembowski and
Piper~\cite[Theorem 4(F) and Corollary 4.9]{dp}. See also
\cite[Section 4.2, statement 10(f)]{dembowski}.

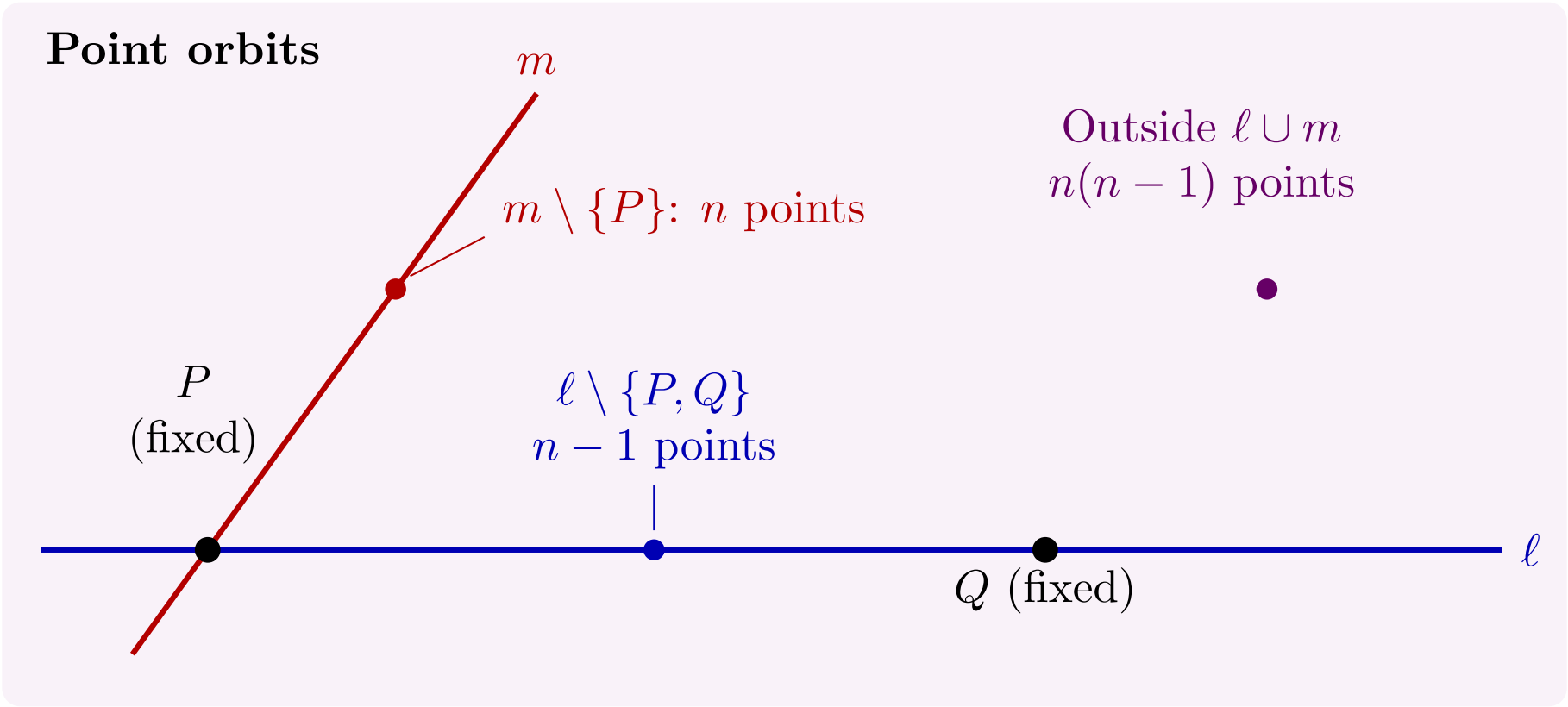
\begin{figure}[H]
\centering
\begin{tikzpicture}[x=1cm,y=1cm,font=\small]
\fill[violet!5,rounded corners] (-1,-1.2) rectangle (11,4.2);
\node[anchor=west] at (-.8,3.85) {\textbf{Point orbits}};
\draw[blue!70!black,very thick] (-.7,0)--(10.5,0) node[right] {$\ell$};
\draw[red!70!black,very thick] (0,-.8)--(3.1,3.5) node[above] {$m$};
\coordinate (P) at (.577,0);
\fill (P) circle (2.8pt) node[above left,align=center,xshift=15pt,yshift=15pt] {$P$\\(fixed)};
\fill (7,0) circle (2.8pt) node[below] {$Q$ (fixed)};
\fill[blue!70!black] (4,0) circle (2.3pt);
\node[blue!70!black,align=center] at (4,1) {$\ell\setminus\{P,Q\}$\\$n-1$ points};
\draw[blue!70!black,thin] (4,.5)--(4,.15);
\fill[red!70!black] (2.018,2) circle (2.3pt);
\node[red!70!black,align=left,anchor=west] at (2.7,2.6) {$m\setminus\{P\}$: $n$ points};
\draw[red!70!black,thin] (2.7,2.4)--(2.13,2.1);
\fill[violet!80!black] (8.7,2) circle (2.3pt);
\node[violet!80!black,align=center] at (8.2,3) {Outside $\ell\cup m$\\$n(n-1)$ points};
\end{tikzpicture}
\caption{The type~(F) configuration and its five point orbits. Each coloured point represents an orbit of the indicated size.}
\label{fig:type-F-point-orbits}
\end{figure}

Ganley~\cite[Theorem 4.3]{ganley} discovered that projective planes
admitting a quasiregular collineation group of type~(F) correspond to what he called
``direct product difference sets''. This correspondence allows the geometric
problem to be studied using group rings, character theory, and the ideal theory of cyclotomic fields,
and is fundamental for this paper.

Let $G=E\times H$ be a finite group where $|E|=n$ and $|H|=n-1$ for some positive integer $n$.
A subset $D\subseteq G$ of cardinality $n-1$ is a \emph{direct product difference set}
(DPDS) of order $n$ if every element of $G\setminus (E \cup H)$
 has exactly one representation $d_1d_2^{-1}$ with $d_1,d_2\in D$,
and no nonidentity element of $E\cup H$ has such a representation.
Equivalently,
\[
DD^{(-1)}=n\,1_G+G-E-H
\]
in the integral group ring $\Z[G]$, with subsets identified with the formal sums over their elements.
We say that $D$ is \emph{abelian} if $G$ is abelian.
Two DPDS $D\subseteq G$ and $D'\subseteq G'$ are \emph{equivalent}
if $D'=g'\alpha(D)$ for some $g'\in G'$ and a group isomorphism
$\alpha:G\to G'$.

By \cite[Theorem 4.3]{ganley}, a DPDS of order $n>2$
in a group $G$ can be used to construct  a projective plane of order $n$ admitting $G$ as 
a quasiregular collineation group of type (F) and vice versa. 
Equivalent DPDS determine isomorphic projective planes 
\cite[Proposition 5.3.1 and p. 134]{pott}.

The classical example of a DPDS is
\begin{equation} \label{classic}
D=\{(x,x):x\in\F_q^*\}\subseteq(\F_q,+)\times(\F_q^*,\cdot),
\end{equation}
where $q$ is a prime power. Ganley attributes this construction to
E.~Spence~\cite[p.~323]{ganley}. Its associated projective plane is the
Desarguesian plane $\operatorname{PG}(2,q)$. See
\cite[Example 5.3.2 and pp.~134--135]{pott} and
\cite[Theorem 1.62]{zhou}.
The main known results on abelian DPDS are the following.

\begin{result}\leavevmode
\begin{enumerate}
\item[(a)] (Pott~\cite{pott94}) Every direct product difference set of
prime order $p$ in an abelian group is equivalent to \eqref{classic}
(with $q=p$).
\item[(b)] (Jungnickel and de Resmini~\cite[Theorem 1]{jr}, via the correspondence of~\cite[Theorem 4.3]{ganley}) If a
 direct product difference set of order $n$ exists in an abelian group
$G=E\times H$, where $|E|=n$ and $|H|=n-1$, then $n=p^a$ for some
prime $p$ and positive integer $a$, and $E$ is elementary abelian.
\end{enumerate}
\end{result}

Jungnickel and de Resmini also stated the conjecture that every
projective plane of order $n$ admitting an abelian collineation group
of order $n(n-1)$ is Desarguesian~\cite[p.~216]{jr}.
In this paper, we prove this conjecture:
\begin{theorem}\label{thm:main}
Every direct product difference set in an abelian group is equivalent
to  \eqref{classic} for some prime power $q$, and hence
every projective plane of order $n$ admitting an abelian collineation group
of order $n(n-1)$ is Desarguesian.
\end{theorem}

A \emph{circular Costas map} is an injective map $\varphi:A\to B$
between finite abelian groups, written additively, with $|B|=|A|+1$,
such that $x\mapsto\varphi(x+h)-\varphi(x)$ is injective for every
$h\in A\setminus\{0\}$. Its graph is a DPDS~\cite{torresrubio},
so our result classifies all such maps.
A \emph{Costas polynomial} is a polynomial $f\in\F_q[X]$ with $f(0)=0$
such that $x\mapsto f(dx)-f(x)$ permutes $\F_q$ for every
$d\in\F_q\setminus\{1\}$.

\begin{corollary}\label{cor:costas}\leavevmode
\begin{enumerate}
\item[(a)] Let $\varphi:A\to B$ be a circular Costas map and let $b$ be the
unique element of $B\setminus\varphi(A)$. Then $q=|B|$ is a prime power,
and there are group isomorphisms $\kappa: A\to(\F_q^*,\cdot)$ and
$\lambda:(\F_q,+)\to B$ such that
\[
\varphi(x)=b+\lambda(\kappa(x))\qquad(x\in A).
\]
In particular, $A$ is cyclic and $B$ is elementary abelian.
\item[(b)] Every Costas polynomial $f\in\F_q[X]$ satisfies
\[
f(X)\equiv L(X^s)\pmod{X^q-X}
\]
for some positive integer $s$ with $\gcd(s,q-1)=1$ and some linearized
permutation polynomial $L\in\F_q[X]$.
\end{enumerate}
\end{corollary}

In the situation of Corollary~\ref{cor:costas}(a), let $g$ be a generator of $A$
and write $\gamma=\phi(g)$. Then $\gamma$ is a primitive element of $\F_q$ and
$\varphi(ig)=b+\lambda(\gamma^i)$ for all integers $i$.
Thus, after translating the omitted value to
zero, circular Costas maps are equivalent to multidimensional Welch maps
$i\mapsto L(\gamma^i)$, where $q=p^a$, $\gamma$ is primitive in $\F_q$, and $L$
is an invertible $\F_p$-linear map.
Corollary~\ref{cor:costas} proves the conjecture of Muratovi\'c-Ribi\'c, Pott,
Thomson and Wang~\cite[Conjecture 4.4]{costas2015} and its extension
to circular Costas maps by Torres and Rubio~\cite{torresrubio}.
For recent work on Costas polynomials, see also~\cite{mrbs}.

\section{Setup and preliminaries}
Throughout the rest of the paper, we assume that $D$ is a DPDS
of prime power order $q$ in an abelian group $G=E\times H$ where $|E|=q$ and $|H|=q-1$.
We write $H$ multiplicatively.
By  the result of Jungnickel and de Resmini, $E$  is elementary abelian,
and we can assume $E=(\F_q,+)$.
The condition $DD^{(-1)}=q\,1_G+G-E-H$ implies that $D$ contains exactly one element
of each coset of $E$ and meets all but one coset of $H$ in exactly one element. 
By translating $D$, if necessary, we can assume that $D\cap H=\emptyset$.
Then
$$
D = \sum_{x\in \F_q^*} x f(x),
$$
where $f: \F_q^*\to H$ is a bijection. 
We extend $f$ to a function on $\F_q$ by adjoining an element $0$ to $H$
and use the convention $f(0)=0$, so that $f$ becomes a bijection $\F_q\to H\cup\{0\}$.

Write  $\zeta_t=e^{2\pi i/t}$ and $q=p^a$ where $p$ is a prime. 
The character values of $D$ lie in the cyclotomic integer ring
\[
\Z[\zeta_{p(q-1)}].
\]
 All congruences below are in this ring. 
We use standard facts on prime ideal factorization and residue class fields of $\Z[\zeta_{p(q-1)}]$,
which can be found in \cite[Chapter~I]{neukirch}, for instance.
Let $\PP$ be a prime ideal  of $\Z[\zeta_{p(q-1)}]$ containing $p$. 
Write $\PP^j\parallel u$ if $u\in\PP^j\setminus\PP^{j+1}$.
Unique factorization of nonzero ideals gives the usual rules for
these exact exponents. 
In particular, $\PP^r\parallel u$ and $uw\in\PP^{r+s}$ imply $w\in\PP^s$
for all $u,w\in\Z[\zeta_{p(q-1)}]$, $r,s\ge 0$.
Moreover, by \cite[Chapter~I, Lemma~10.1 and Proposition~10.3]{neukirch}, we have
\begin{equation}\label{eq:arithmetic}
\begin{gathered}
\Z[\zeta_{p(q-1)}]/\PP\cong\F_q,\\
\PP\parallel(\zeta_p-1),\qquad
\PP^{p-1}\parallel p,\qquad
\PP^{a(p-1)}\parallel q.
\end{gathered}
\end{equation}

We use the character orthogonality relations. In particular, for every
complex character $\chi$ of $H$,
\[
\sum_{h\in H}\chi(h)=
\begin{cases}
|H|,&\text{if $\chi$ is the trivial character},\\
0,&\text{if $\chi$ is nontrivial}.
\end{cases}
\]
Recall that  $f(0)=0$ for the bijection $f:\F_q\to H\cup \{0\}$.
Using the convention $\chi(0)=0$ for a character $\chi$ of $H$, we get
\begin{equation}\label{eq:character-sum-zero}
\sum_{x\in\F_q}\chi(f(x))=\sum_{h\in H}\chi(h)=0
\qquad\text{for nontrivial }\chi.
\end{equation}
\section{The character values of $D$}

Recall that the DPDS $D$ is uniquely determined by the bijection
$f: \F_q \to H\cup\{0\}$ (with the convention $f(0)=0$).
In this section, we study the character values of $D$ by  using Fourier inversion over $(\F_q,+)$.

Write $\Tr=\Tr_{\F_q/\F_p}$. A value of the trace in an exponent is
represented by any integer in its residue class modulo $p$.
Fix a nontrivial complex character $\chi$ of $H$.
The additive characters of $\F_q$ are precisely
$x\mapsto\zeta_p^{\Tr(xy)}$, with $y\in\F_q$.
Thus, for each $y$, the map
$xh\mapsto\zeta_p^{\Tr(xy)}\chi(h)$, $x\in\F_q$, $h\in H$, is a character of
$(\F_q,+)\times H$. Its value on $D$ is
\[
a_y=\sum_{x\in\F_q}\zeta_p^{\Tr(xy)}\chi(f(x))
\qquad(y\in\F_q).
\]
The term at $x=0$ vanishes by our convention $f(0)=0$ and $\chi(0)=0$.
Hence the $a_y$'s are
exactly the character values of $D$ obtained by fixing $\chi$ and varying
the additive character.  Note that  $a_y\in \Z[\zeta_{p(q-1)}]$ for all $y$ and
$a_0=0$ by \eqref{eq:character-sum-zero}.
For $y\ne0$, both factors of the product character are nontrivial.
Its values on the group-ring sums $G$, $E$, and $H$ therefore vanish.
Applying it to $DD^{(-1)}=q\,1_G+G-E-H$ gives
\begin{equation}\label{eq:norm}
a_y\overline{a_y}=q\qquad(y\ne0).
\end{equation}

For each element of $\F_p$, fix an integer representative modulo $p$,
with $0$ represented by $0$.
We use these representatives when trace values occur as coefficients
in sums of algebraic integers. Changing a representative
changes the sum by $p$ times an algebraic integer, which belongs to
$\PP$. Thus the residue class of the sum is independent of these choices.
By \eqref{eq:character-sum-zero} and since $f$ is a bijection, we have 
\begin{equation}\label{e1}
a_y=\sum_{x\in\F_q}\bigl(\zeta_p^{\Tr(xy)}-1\bigr)\chi(f(x)).
\end{equation}
For $1\le j<p$,
\[
\frac{\zeta_p^j-1}{\zeta_p-1}
=1+\zeta_p+\cdots+\zeta_p^{j-1}\equiv j\pmod{\PP},
\]
and for $j=0$ the quotient is zero.
Thus every summand in \eqref{e1} is divisible by $\zeta_p-1$.
In particular, $a_y/(\zeta_p-1)$ is an algebraic integer, and
\begin{equation}\label{eq:quotient-residue}
\frac{a_y}{\zeta_p-1}
\equiv\sum_{x\in\F_q}\Tr(xy)\chi(f(x))\pmod{\PP}.
\end{equation}

For any function $u:\F_q\to\mathbb C$,  Fourier inversion on $(\F_q,+)$ \cite[Chapter VI, Lemma 3.5, p.~314]{bjl} gives 
\begin{equation}\label{eq:fourier-inversion-additive}
u(z)=\frac1q\sum_{y\in\F_q}
\left(\sum_{x\in\F_q}u(x)\zeta_p^{-\Tr(xy)}\right)
\zeta_p^{\Tr(yz)}
\qquad(z\in\F_q).
\end{equation}
 Note that \eqref{eq:fourier-inversion-additive}
is an identity over $\mathbb C$, before any reduction
modulo $\PP$.
For $u(x)=\overline{\chi(f(x))}$, the inner sum on the right hand side of \eqref{eq:fourier-inversion-additive} is $\overline{a_y}$.
Multiplying \eqref{eq:fourier-inversion-additive} by $q$ gives
\begin{equation}\label{eq:fourier}
q\,\overline{\chi(f(z))}
=\sum_{y\in\F_q}\zeta_p^{\Tr(yz)}\overline{a_y}.
\end{equation}
At $z=0$, the left side is zero, so
$\sum_{y\in\F_q}\overline{a_y}=0$.
Subtracting this identity from \eqref{eq:fourier} yields
\begin{equation}\label{eq:difference}
q\,\overline{\chi(f(z))}=\sum_{y\in\F_q}
\bigl(\zeta_p^{\Tr(yz)}-1\bigr)\overline{a_y}.
\end{equation}

\section{The normalized character values}
We have shown in the previous section that all character values $a_y$
are divisible by $\zeta_p-1$ and thus we can consider the ``normalized'' values $a_y/(\zeta_p-1)$. 
In this section, we study these values -- this time with Fourier inversion on $H$.
For any function $u:H\to\mathbb C$, Fourier inversion \cite[Chapter VI, Lemma 3.5, p.~314]{bjl}  takes the form
\begin{equation}\label{eq:fourier-inversion-H}
u(h)=\frac{1}{|H|}\sum_{\chi}
\left(\sum_{g\in H}u(g)\chi(g)\right)\chi(h^{-1})
\qquad(h\in H),
\end{equation}
where the outer sum runs over all complex characters of $H$.
The inner sum is the Fourier coefficient at $\chi$ in this convention.
Indeed, after interchanging the two sums, character orthogonality gives
$\sum_\chi\chi(gh^{-1})=|H|$ if $g=h$, and zero otherwise.
In particular, if every Fourier coefficient vanishes, then $u$ is
identically zero.

For a fixed nontrivial character $\chi$ of $H$, define
\[
F_\chi:\F_q\to\Z[\zeta_{p(q-1)}],\ y \mapsto \frac{a_y}{\zeta_p-1}.
\]
The preceding section shows that $F_\chi(y)$ indeed is  an algebraic integer for all $y$.

\begin{lemma}\label{lem:first}
We have
\begin{equation}\label{eq:F}
F_\chi(y)\equiv\sum_{x\in\F_q}\Tr(xy)\chi(f(x))
\pmod{\PP}
\end{equation}
for all $y$.
The map 
$$
 \F_q \to \Z[\zeta_{p(q-1)}]/\PP,\  y\mapsto F_\chi(y)+\PP
 $$ is $\F_p$-linear.
If $q>2$, then there exists a nontrivial $\chi$ for which $F_\chi(1)\notin\PP$.
\end{lemma}

\begin{proof}
Note that  \eqref{eq:F} follows from \eqref{eq:quotient-residue} and the definition of $F_\chi$.
Linearity of $y\mapsto F_\chi(y)+\PP$ follows from linearity of the trace
(even though we chose integer representatives for trace values, 
the linearity over $\F_p$ suffices here since $p\in \PP$).

To prove the last statement, consider the function
$h\mapsto\Tr(f^{-1}(h))$ from $H$ to $\Z$
(recall our convention to represent trace values by integers).
This function is not identically zero ${\rm mod}\ p$ because the trace is a nonzero linear form.
Note that
\begin{equation} \label{fourier_trace}
|H| \Tr(f^{-1}(h)) =\sum_{\chi}
\left(\sum_{g\in H}\Tr(f^{-1}(g))\chi(g)\right)\chi(h^{-1})
\qquad(h\in H)
\end{equation}
by \eqref{eq:fourier-inversion-H}.
Moreover, we have
\begin{equation} \label{fourier_trace1}
\sum_{g\in H}\Tr(f^{-1}(g))
=\sum_{x \in \F_q}\Tr(x)
\equiv 0\pmod p,
\end{equation}
since $\Tr(0)=0$, $\sum_{x \in \F_q}x=0$ in $\F_q$ for $q>2$, and the trace
is $\F_p$-linear.
On the other hand,
\begin{equation} \label{fourier_trace2}
\sum_{g\in H}\Tr(f^{-1}(g))\chi(g)
= \sum_{x\in \F_q}\Tr(x)\chi(f(x))
\equiv F_\chi(1)\pmod{\PP}
\end{equation}
for every nontrivial character $\chi$ of $H$ by \eqref{eq:F}.
Note that $|H|=q-1 \not\in \PP$.
Thus, if $F_\chi(1)\in\PP$ for every nontrivial $\chi$, 
then  $\Tr(f^{-1}(h)) \in \PP$ for all $h$  by \eqref{fourier_trace}--\eqref{fourier_trace2}.
As $\PP\cap\Z=p\Z$ and $\Tr(f^{-1}(h)) \in \Z$, this implies
$\Tr(f^{-1}(h)) \equiv 0\pmod p$ for all $h$,
a contradiction.
\end{proof}

\section{Kernel dichotomy}
We now show that, for every nontrivial character $\chi$ of $H$,
the kernel of $y\mapsto F_\chi(y)+\PP$ is either $\{0\}$ or $\F_q$.

\begin{lemma}\label{lem:kernel}
For every nontrivial character $\chi$ of $H$, the map
\[
\F_q\to\Z[\zeta_{p(q-1)}]/\PP,\ y\mapsto F_\chi(y)+\PP
\]
is either zero or bijective.
\end{lemma}

\begin{proof}
Suppose the map $y\mapsto F_\chi(y)+\PP$ is nonzero, and let
\[
V=\{y\in\F_q:F_\chi(y)\in\PP\}.
\]
Then $V$ is a proper $\F_p$-subspace of $\F_q$.
For $y\notin V$, the identity $a_y=(\zeta_p-1)F_\chi(y)$
and \eqref{eq:norm} imply
\begin{equation}\label{eq:orders}
\PP\parallel a_y,\qquad \PP^{a(p-1)-1}\parallel\overline{a_y}.
\end{equation}
If $y,y'\notin V$ belong to the same coset of $V$, then
$F_\chi(y)\equiv F_\chi(y')\pmod{\PP}$,
so
\[
a_{y'}-a_y\in(\zeta_p-1)\PP\subseteq\PP^2.
\]
Moreover, by \eqref{eq:norm},
\[
a_ya_{y'}(\overline{a_y}-\overline{a_{y'}})
=q(a_{y'}-a_y).
\]
The right side belongs to $\PP^{a(p-1)+2}$, while $\PP^2\parallel a_ya_{y'}$. 
We conclude
\begin{equation}\label{eq:coset}
\overline{a_y}\equiv\overline{a_{y'}}
\pmod{\PP^{a(p-1)}}.
\end{equation}

Suppose that $V\ne\{0\}$.
Since $V$ is a proper subspace of $\F_q$, there is  a nonzero $t\in\F_q$ such that
\[
\Tr(tv)=0\text{ in }\F_p\qquad\text{for every }v\in V.
\]
Note that 
\begin{equation} \label{dicho1}
\zeta_p^{\Tr(ty)}-1=0 \text{ for } y\in V.
\end{equation}
For $y$ in a  coset $y_0+V \neq V$, the factor $\zeta_p^{\Tr(ty)}-1$ is
constant, equal to $\zeta_p^{\Tr(ty_0)}-1$, and belongs to $\PP$.
By \eqref{eq:coset}, we have
\[
\sum_{y\in y_0+V}\overline{a_y}
\equiv |V|\,\overline{a_{y_0}}\pmod{\PP^{a(p-1)}}.
\]
Since $V\ne\{0\}$, its order is divisible by $p$, so $|V|\in\PP$.
Together with $\overline{a_{y_0}}\in\PP^{a(p-1)-1}$, this shows that
each coset sum belongs to $\PP^{a(p-1)}$. Consequently
\begin{equation} \label{dicho2}
\sum_{y\in y_0+V}\bigl(\zeta_p^{\Tr(ty)}-1\bigr)\overline{a_y}
\in\PP^{a(p-1)+1}.
\end{equation}
However,
$$
q\,\overline{\chi(f(t))}
=\sum_{y\in\F_q} \bigl(\zeta_p^{\Tr(yt)}-1\bigr)\overline{a_y}
$$
by \eqref{eq:difference}.
Thus,  $q\,\overline{\chi(f(t))}\in\PP^{a(p-1)+1}$ by \eqref{dicho1} and \eqref{dicho2}.
This is impossible: $t\ne0$, so $\overline{\chi(f(t))}$ is a root
of unity and a unit, and $\PP^{a(p-1)}\parallel q\,\overline{\chi(f(t))}$.

Therefore $V=\{0\}$, so the map $y\mapsto F_\chi(y)+\PP$ is injective.
Its domain and codomain both have exactly $q$ elements, hence it is bijective.
\end{proof}

\section{Reconstruction and classification} \label{recon}
We fix an identification $\Z[\zeta_{p(q-1)}]/\PP\cong\F_q$.

\begin{theorem}\label{thm:normalized}
Let $p$ be a prime, $a\ge1$, $q=p^a$, and let $H$ be an abelian group
of order $q-1$. Suppose that $f:\F_q^*\to H$ is a bijection and
\[
D=\sum_{x\in\F_q^*}xf(x)
\]
is a DPDS in $(\F_q,+)\times H$. Then there are an additive automorphism
$\alpha:\F_q\to\F_q$ and a group isomorphism
$\phi:H\to\F_q^*$ satisfying
\[
\phi(f(x))=\alpha(x)\qquad(x\in\F_q^*).
\]
Consequently $H$ is cyclic, and the group isomorphism
$(\F_q,+)\times H\to(\F_q,+)\times(\F_q^*,\cdot)$ induced by $\alpha$ and $\phi$ maps
$D$ to \eqref{classic}.
\end{theorem}

\begin{proof}
The case $q=2$ is trivial. 
For $q>2$, choose a nontrivial
character $\chi$ for which $F_\chi(1)\notin\PP$, as in
Lemma \ref{lem:first}. Lemma \ref{lem:kernel} gives
$F_\chi(y)\notin\PP$ for every $y\ne0$.
Thus each algebraic integer $F_\chi(y)$, for $y\ne0$, is nonzero modulo $\PP$. 
Recall that $F_\chi(y)=a_y/(\zeta_p-1)$.
Multiply \eqref{eq:difference} by
$\prod_{y\in\F_q^*}F_\chi(y)$, use $a_y\overline{a_y}=q$, and
cancel $q$. This gives the following identity of algebraic integers:
\begin{equation}\label{eq:cleared}
\overline{\chi(f(z))}
\prod_{y\in\F_q^*}F_\chi(y)
=\sum_{y\in\F_q}\frac{\zeta_p^{\Tr(yz)}-1}{\zeta_p-1}
\prod_{\substack{w\in\F_q^*\\w\ne y}}F_\chi(w).
\end{equation}
Reducing modulo $\PP$ and dividing by the nonzero element
$\left(\prod_{y\in\F_q^*}F_\chi(y)\right)+\PP$ of the residue field gives
\begin{equation}\label{eq:alpha}
\alpha(z):=\overline{\chi(f(z))}+\PP
=\sum_{y\in\F_q^*}\frac{\Tr(yz)}{F_\chi(y)+\PP}
\qquad(z\in\F_q).
\end{equation}
Recall that $\Z[\zeta_{p(q-1)}]/\PP\cong \F_q$ and note that
\eqref{eq:alpha} is an identity in this residue field. 
The right side is $\F_p$-linear in $z$.
Also $\alpha(0)=0$ and $\alpha(z)\ne0$ for $z\ne0$, since
$\overline{\chi(f(z))}$ is a root of unity.
Thus $\alpha:\F_q\to\F_q$ is an additive automorphism.

Define
\[
\phi:H\to\F_q^*,\
\phi(h)=\overline{\chi(h)}+\PP.
\]
This is a group homomorphism, and \eqref{eq:alpha} says
$\phi(f(z))=\alpha(z)$ for $z\ne0$.
Since $f$ is bijective and $\alpha$ maps $\F_q^*$ onto itself,
$\phi$ is a multiplicative isomorphism. 
These maps have  the asserted properties.
\end{proof}

\begin{proof}[Proof of Theorem~\ref{thm:main}]
Let $D$ have order $n$ in the abelian group $G=E\times H$,
where $|E|=n$ and $|H|=n-1$.
The result of Jungnickel and de Resmini \cite{jr}, applied through
the DPDS correspondence, gives $n=p^a$ and makes $E$ elementary
abelian. Identify $E$ with $(\F_q,+)$, where $q=p^a$.
After translating $D$, we may assume that
$D=\sum_{x\in\F_q^*}xf(x)$ for a bijection $f:\F_q^*\to H$.
Theorem~\ref{thm:normalized} shows that $D$ is equivalent to \eqref{classic}.

Finally, the unique projective plane of order $2$ is $\operatorname{PG}(2,2)$.
For $n>2$, let $G$ be an abelian collineation group of order $n(n-1)$ of a
projective plane of order $n$. Then $G$ is quasiregular and thus has an orbit
decomposition of one of the types (A)--(F) by \cite[Corollary 4.9]{dp}.
Among these types, only type~(F) has $|G|=n(n-1)$.
Hence the plane corresponds to an abelian DPDS
by \cite[Theorem 4.3]{ganley}. Equivalent DPDS determine isomorphic
projective planes, and the classical DPDS determines
$\operatorname{PG}(2,q)$~\cite[Section 5.3]{pott}.
Hence the plane is Desarguesian.
\end{proof}

\begin{proof}[Proof of Corollary~\ref{cor:costas}]
(a) Since $\varphi$ is injective and $|B|=|A|+1$, there is exactly one
element $b$ of $B$ that is not contained in $\varphi(A)$.
The graph of $\varphi$ is a DPDS of order $|B|$ in $B\times A$~\cite{torresrubio}.
By the result of Jungnickel and de Resmini, $q=|B|$ is a prime power and
$B$ is elementary abelian. Fix a group isomorphism $\beta:B\to(\F_q,+)$.
Then $x\mapsto\beta(\varphi(x)-b)$ is a bijection $A\to\F_q^*$.
Let $f:\F_q^*\to A$ be its inverse.
Translating the graph of $\varphi$ by $-b$ and applying $\beta$ to its
$B$-components, we obtain the DPDS $\sum_{z\in\F_q^*}zf(z)$ in $(\F_q,+)\times A$.
By Theorem~\ref{thm:normalized}, applied with $H=A$, there are an additive
automorphism $\alpha$ of $\F_q$ and a group isomorphism $\kappa: A\to\F_q^*$ with
$\kappa(f(z))=\alpha(z)$ for all $z\in\F_q^*$.
Setting $z=\beta(\varphi(x)-b)$, we get $\kappa(x)=\alpha(\beta(\varphi(x)-b))$ and thus
$\varphi(x)=b+\lambda(\kappa(x))$ for all $x\in A$, where $\lambda=\beta^{-1}\circ\alpha^{-1}$.
As $\F_q^*$ is cyclic, so is $A$.

(b) Let $f\in\F_q[X]$ be a Costas polynomial. The condition for $d=0$ shows that
$f$ permutes $\F_q$. As $f(0)=0$, the restriction $\varphi$ of $f$ to $\F_q^*$
is an injective map from $A=(\F_q^*,\cdot)$ to $B=(\F_q,+)$, and $0$ is the
unique element of $B\setminus\varphi(A)$. For every $d\in\F_q^*\setminus\{1\}$,
the map $x\mapsto\varphi(dx)-\varphi(x)$ is injective on $\F_q^*$.
Hence $\varphi$ is a circular Costas map (with $A$ written multiplicatively).
By part (a), there are an automorphism $\kappa$ of $(\F_q^*,\cdot)$ and
an automorphism $\lambda$ of $(\F_q,+)$ with $f(x)=\lambda(\kappa(x))$ for all $x\in\F_q^*$.
There is a positive integer $s$ with $\gcd(s,q-1)=1$ and $\kappa(x)=x^s$ for all
$x\in\F_q^*$, and $\lambda$ is induced by a linearized permutation polynomial $L\in\F_q[X]$.
Since $f(0)=0=L(0)$, we conclude $f(x)=L(x^s)$ for all $x\in\F_q$, which is
equivalent to $f(X)\equiv L(X^s)\pmod{X^q-X}$.
\end{proof}

\section*{AI statement}
The key approach used in this paper, 
employing the residue class field to reconstruct the classical DPDS in $\F_q \cong \Z[\zeta_{p(q-1)}]/\PP$
and the kernel dichotomy,
is solely due to the authors who started this work in characteristic $2$. 
 Claude Fable 5.2 made substantial contributions to the reconstruction (Section \ref{recon}) in characteristic $2$,
 and GPT-6 Astra to the generalization to arbitrary characteristic.
All proofs have been verified by the authors.

\end{document}